\input amstex
\input amsppt.sty

\define\wt{\widetilde}
\define\p{\Cal P}
\define\q{\overline{\Cal P}}
\define\const{\operatorname{const}}
\define\tht{\thetag}

\NoBlackBoxes
\TagsOnRight

%\magnification 1200
\topmatter 
\title
Duality of orthogonal polynomials on a finite set
\endtitle
\author
Alexei Borodin
\endauthor

%\date Preliminary version, January 14, 2000
%\enddate

\address
Department of Mathematics, The University of
Pennsylvania, Philadelphia, PA 19104-6395, U.S.A.  
E-mail address:
{\tt borodine\@math.upenn.edu}
\endaddress

\abstract
We prove a certain duality relation for orthogonal polynomials defined on a
finite set. The result is used in a direct proof of the equivalence of two
different ways of computing the correlation functions of a discrete orthogonal
polynomial ensemble.
\endabstract
\endtopmatter

\document

\head Introduction
\endhead

This note is about a certain duality of orthogonal polynomials defined
on a finite set. If the weights of two systems of orthogonal
polynomials are related in a certain way, then the values of the $n$th
polynomial of the first system at the points of the set equal, up to a simple
factor, the corresponding values of the $(M-n)$th polynomial of the second system,
where $M$ is the cardinality of the underlying finite set. 

We formulate the exact result and prove it in \S1.

In \S2 we explain the motivation which led to the result. We compare two
different ways to compute probabilistic quantities called
{\it correlation functions} in a certain model. The model is a
discrete analog of the {\it orthogonal polynomial ensembles} which appeared
for the first time in the random matrix theory, see, e.g.,
\cite{Dy}, \cite{Ga}, \cite{GM}, \cite{Me}, \cite{NW}.  Discrete orthogonal polynomial
ensembles  were discussed in \cite{BO1}, \cite{BO2}
\cite{BO3}, \cite{J1}-\cite{J3}.

The results of the two computations must be
equal, but this is not at all obvious from the explicit formulas. Our duality
relation provides a proof of the equivalence of the two
resulting expressions.     

In \S3 we consider 2 examples when the orthogonal polynomials are
classical (Krawtchouk and Hahn polynomials). In these cases the
duality provides relations between similar polynomials with different sets of
parameters. The relations are also easily verified using known explicit
formulas for the polynomials. 

I am very grateful to Grigori Olshanski for numerous discussions. I also want
to thank Tom Koornwinder for providing me with his computation regarding the
Hahn polynomials, see \S3.  

\head 1. Duality
\endhead

\proclaim{Theorem 1} Let 
$$
X=\{x_0,x_1,\dots,x_M\}\subset \Bbb R
$$ 
be a finite set of
distinct points on the real line, $u(x)$ and $v(x)$ be two positive functions
on $X$ such that  
$$ 
u(x_k)v(x_k)=\frac 1{\prod_{i\ne k}(x_k-x_i)^2}\,,
\quad k=0,1,\dots,M,  
\tag 1
$$
and  $P_0,P_1,\dots,P_{M}$ and $Q_0,Q_1,\dots,Q_M$ be the systems of
orthogonal polynomials on $X$ with respect to the weights $u(x)$ and
$v(x)$, respectively, 
$$ 
\gathered
\deg P_i=\deg Q_i=i,\\
\sum_{k=0}^M P_i(x_k)P_j(x_k)u(x_k)=\delta_{ij}p_i,\quad 
\sum_{k=0}^M Q_i(x_k)Q_j(x_k)v(x_k)=\delta_{ij}q_i,   
\\ P_i=a_ix^i+\text{\rm lower terms},\quad
Q_i=b_ix^i+\text{\rm lower terms}. \endgathered
$$

Assume that  the polynomials are normalized so that $b_i=p_{M-i}/a_{M-i}$ for
all $i=0,1,\dots,M$.
Then
$$
\gathered
P_i(x)\sqrt{u(x)}=\epsilon(x)Q_{M-i}(x)\sqrt{v(x)},\quad x\in X,\\
a_ib_{M-i}=p_i=q_{M-i}, \quad i=0,1,\dots,M,
\endgathered
$$
where
$$
\epsilon(x_k)=\operatorname{sgn}\prod_{i\ne k} (x_k-x_i),\quad k=0,1,\dots,M.
$$
\endproclaim
\demo{Proof}
Let us start  with one system of polynomials, say, $\{P_i\}$, and define a
sequence of functions $\{\wt Q_i\}$ on $X$ by the equalities 
$$
\wt
Q_i(x_k)=\epsilon(x_k)P_{M-i}(x_k)\sqrt{\frac{u(x_k)}{v(x_k)}}=\prod_{i\ne k}
(x_k-x_i)\cdot P_{M-i}(x_k){u(x_k)}\,. 
$$

Then 
$$
 \sum_{k=0}^M \wt Q_i(x_k)\wt Q_j(x_k)v(x_k)=\sum_{k=0}^M
P_{M-i}(x_k)P_{M-j}(x_k) u(x_k)=\delta_{ij}p_{M-i}, 
$$
so the functions $\{\wt Q_i\}_{i=0}^M$ are pairwise orthogonal with respect
to the weight $v(x)$, and 
$
q_i=\Vert \wt Q_i \Vert^2_v=p_{M-i}.
$

 Consider the interpolation polynomial $Q_i(x)$ of degree $M$ such that
$Q_i(x)=\wt Q_i(x)$ for all $x\in X$. We have (the hat means that the
corresponding factor is omitted) 
$$ 
\multline
 Q_i(x)=\sum_{m=0}^M \wt Q_i(x_m)\,\frac{(x-x_0)\cdots
\widehat{(x-x_m)}\cdots (x-x_M)}{(x_m-x_0)\cdots \widehat{(x_m-x_m)}\cdots
(x_m-x_M)}\\ =\sum_{m=0}^M P_{M-i}(x_m)u(x_m)\cdot (x-x_0)\cdots
\widehat{(x-x_m)}\cdots (x-x_M).  
\endmultline 
$$
The coefficient of $x^n$ of such polynomial equals
$$
(-1)^{M-n}\sum_{m=0}^M P_{M-i}(x_m)u(x_m)e_{M-n}(x_0,\dots,\widehat{x_m},\dots,x_M)
$$
where $e_s$ are the elementary symmetric functions:
$$
e_s(y_0,y_1,\dots)=\sum_{0\le i_1<i_2<\dots<i_s} y_{i_1}y_{i_2}\cdots y_{i_s}.
$$
Denote $e_s(x_0,\dots,x_M)$ by $E_s$. Note that $E_0=1$ by
definition. An application of the inclusion--exclusion principle shows that 
$$
e_s(x_0,\dots,\widehat{x_m},\dots,x_M)=E_s-x_mE_{s-1}+x_m^2E_{s-2}-\dots+(-1)^sx_m^s.
$$
Then the coefficient of $x^n$ in $Q_i(x)$ equals
$$
\multline
(-1)^{M-n}\sum_{m=0}^M
 P_{M-i}(x_m)u(x_m)\left(E_{M-n}-x_mE_{M-n-1}+\dots+(-1)^{M-n}x_m^{M-n}\right)\\=
(-1)^{M-n}E_{M-n}\langle P_{M-i}, 1\rangle +(-1)^{M-n-1}E_{M-n-1}\langle
P_{M-i}, x\rangle+\dots+\langle P_{M-i}, x^{M-n}\rangle. 
\endmultline 
$$
But the orthogonality of $P_j$'s implies that
$
\langle P_{M-i}, x^r\rangle=0$ for  $r<M-i$, and
$$
\langle P_{M-i}, x^{M-i}\rangle=\frac{\Vert P_{M-i}\Vert ^2}{a_{M-i}}
=\frac{p_{M-i}}{a_{M-i}}\,. 
$$
 This immediately implies that $Q_i$ is a polynomial of degree $i$ with the
leading coefficient $b_i={p_{M-i}}/{a_{M-i}}$. \qed 
\enddemo

\head 2. Probabilistic interpretation
\endhead

Recall that $X=\{x_0,\dots, x_M\}$ is a finite subset of the real line. 

For any $m=1,\dots, M$, denote by $X^{(m)}$ the set of all
subsets of $X$ with $m$ points:
$$
X^{(m)}=\left\{\{x_{i_1},\dots,x_{i_m}\}\,|\, 0\le i_1<\dots<i_m\le M\right\}.
$$

For any positive function $w(x)$ on $X$ denote by $\p_w^{(m)}$ the probability
measure on $X^{(m)}$ defined by the formula:
$$
\p_w^{(m)}\{x_{i_1},\dots,x_{i_m}\}=\const\prod_{1\le k<l\le
m}(x_{i_k}-x_{i_l})^2\cdot \prod_{k=1}^m w(x_{i_k}). 
$$

Also denote by $\q_w^{(m)}$ the probability
measure on $X^{(m)}$ defined by the relation:
$$
\q_w^{(m)}(A)=\p_w^{(m)}(X\setminus A),\quad A\in X^{(m)}. 
$$

The next claim was essentially proved in \cite{BO3}.

\proclaim{Proposition 2} Let $u(x)$ and $v(x)$ be two positive functions on
$X$ satisfying \tht{1}. Then $\p_u^{(m)}=\q_v^{(M-m+1)}$ for any $m=1,\dots,M$.
\endproclaim 

\demo{Proof}
For arbitrary finite sets $B$ and $C$ we will abbreviate
$$
\Pi(B)=\pm\prod_{\Sb x,y\in B\\ x\ne y\endSb}(x-y),\quad
\Pi(B,C)=\prod_{x\in B,\,y\in C}(x-y).
$$
The sign of $\Pi(B)$ is inessential. 

Take $A=\{x_{i_1},\dots, x_{i_m}\}\in X^{(m)}$.
We have 
$$
\p_u^{(m)}(A)=\const\prod_{1\le k<l\le
m}(x_{i_k}-x_{i_l})^2\cdot \prod_{k=1}^m u(x_{i_k})=
\const\cdot\Pi^2(A)\cdot \prod_{x\in A}
u(x). 
$$

Further,
$$
\Pi(A)=\pm\Pi(X\setminus A)\cdot{\Pi^2(A)
\Pi(A,X\setminus A)}\cdot\frac 1{  
\Pi(A)\Pi(X\setminus A)\Pi(A,X\setminus A)}\,.
$$
But $\Pi(A)\Pi(X\setminus A)\Pi(A,X\setminus A)=\Pi(X)=\const$, and
$$
\Pi^2(A)\Pi(A,X\setminus A)=\pm \prod_{x\in A}\left(\prod_{\Sb y\in X\\ y\ne
x\endSb} (y-x)\right). 
$$
Hence, using \tht{1}, we get
$$
\gathered
\Pi^2(A)\cdot \prod_{x\in A} u(x)=\const \cdot\Pi^2(X\setminus
A)\left(\prod_{x\in A}v(x)\right)^{-1}\\=\const \cdot\Pi^2(X\setminus
A)\cdot\frac{\prod_{x\in X\setminus A}v(x)}{\prod_{x\in X}v(x)}
=\const' \cdot\Pi^2(X\setminus A)\cdot\prod_{x\in X\setminus A}v(x), 
\endgathered
$$ 
where $\const'=\const\cdot \left(\prod_{x\in X}v(x)\right)^{-1}$.
Thus, $\p_u^{(m)}$ and $\q_v^{(M-m+1)}$ differ by a multiplicative constant.
Since both $\p_u^{(m)}$ and $\q_v^{(M-m+1)}$ are probability measures, they must
coincide. \qed
\enddemo

Let $\mu$ be an arbitrary probability measure on the set of all subsets of
$X$. Note that any probability measure on $X^{(m)}$ can be trivially extended
to a measure on the set of all subsets of $X$.

For any $n=1,2,\dots,M$, we define the $n$th correlation function of $\mu$
$$
\rho_n(\,\cdot\,|\,\mu): X^{(n)}\to \Bbb R_{\ge 0}
$$
by the formula
$$
\rho_n(A\,|\, \mu)=\sum_{B\supset A}\mu(B).
$$ 
In other words, $\rho_n(A\,|\,\mu)$ is the probability (with respect to
$\mu$) that the random set $B$ contains a fixed set $A\in X^{(n)}$.

Below we use the notation of Theorem 1 for the orthogonal polynomials
associated with the weights $u(x)$ and $v(x)$. 

\proclaim{Proposition 3} For any $m=1,\dots,M$, the correlation functions
of $\p^{(m)}_u$ have the form
$$
\rho_n(\{x_{i_1},\dots,x_{i_n}\}\,|\,
\p^{(m)}_u)=\det\left[K_u^{(m)}\left(x_{i_k},x_{i_l}\right)\right]_{k,l=1,\dots,
n}, 
$$
where
$$
K_u^{(m)}\left(x,y\right)=\sqrt{u(x)u(y)}\,\sum_{i=0}^{m-1}
\frac{P_i(x)P_i(y)}{p_i}\,.
$$
\endproclaim
\demo{Proof} A standard argument from the random matrix theory, see, e.g,
\cite{Dy}, \cite{Me, 5.2}. \qed
\enddemo

Note that if $n>m$ then the $n$th correlation function of $\p^{(m)}_u$
vanishes identically. Indeed, all sets with more than $m$ points have measure
zero with respect to $\p^{(m)}_u$. Another way to see the vanishing is to
observe that the matrix $\Vert K_u^{(m)}(x_i,x_j)\Vert_{i,j=0,\dots,M}$ 
has rank $m$. Thus, its $n\times n$ minors expressing
$\rho_n(\,\cdot\,|\,\p^{(m)}_u)$ must vanish if $n>m$. 

Similarly, for any $m=1,\dots,M$, the correlation functions of $\p^{(m)}_v$
have the form
$$
\rho_n(\{x_{i_1},\dots,x_{i_n}\}\,|\,
\p^{(m)}_v)=\det\left[K_v^{(m)}\left(x_{i_k},x_{i_l}\right)\right]_{k,l=1,\dots,
n}, 
$$
where
$$
K_v^{(m)}\left(x,y\right)=\sqrt{v(x)v(y)}\,\sum_{i=0}^{m-1}
\frac{Q_i(x)Q_i(y)}{q_i}\,.
$$

The determinantal formulas for the correlation functions above imply that
$\p_u^{(m)}$ and $\p_v^{(m)}$ belong to the class of {\it determinantal point
processes}, see \cite{Ma}, \cite{DVJ, 5.4}, \cite{BOO, Appendix}, \cite{So}
for a general discussion of such processes.

\proclaim{Proposition 4} For any $m=1,\dots,M$, the correlation
functions of $\q^{(m)}_u$ have the form
$$
\rho_n(\{x_{i_1},\dots,x_{i_n}\}\,|\,
\q^{(m)}_u)=\det\left[
\overline{K}_u^{(m)}\left(x_{i_k},x_{i_l}\right)\right]_{k,l=1,\dots,
n}, 
$$
where
$$
\overline{K}_u^{(m)}\left(x,y\right)=\delta_{xy}-K_u^{(m)}\left(x,y\right).
$$
Here $\delta_{xy}$ is the Kronecker delta.
\endproclaim
\demo{Proof} By the definition of $\q_u^{(m)}$,
we have
$$
\rho_n(A\,|\,\q^{(m)}_u)=\sum_{B\supset A}\p_u^{(m)}(X\setminus B)=
\sum_{\Sb C\subset X\\ C\cap A=\varnothing\endSb}\p_u^{(m)}(C).
$$
The inclusion-exclusion principle, see, e.g., \cite{Ha, 2.1}, gives 
$$
\sum_{\Sb C\subset X\\ C\cap A=\varnothing\endSb}\p_u^{(m)}(C)=
\sum_{D\subset A}(-1)^{|D|}\rho_{|D|}(D\,|\,\p_u^{(m)}).
$$
By Proposition 3, the expression on the right-hand side is equal to the
alternating sum of all diagonal minors of the matrix 
$\Vert K_u^{(m)}(x,y)\Vert_{x,y\in A}$. By linear algebra, this is equal to 
$\det[\delta_{xy}- K_u^{(m)}(x,y)]_{x,y\in A}$. \qed
\enddemo

Similarly, for any $m=1,\dots,M$, the correlation functions of $\q^{(m)}_v$
have the form
$$
\rho_n(\{x_{i_1},\dots,x_{i_n}\}\,|\,
\q^{(m)}_v)=\det\left[\overline{K}_v^{(m)}
\left(x_{i_k},x_{i_l}\right)\right]_{k,l=1,\dots,
n}, 
$$
where
$$
\overline{K}_v^{(m)}\left(x,y\right)=\delta_{xy}-K_v^{(m)}\left(x,y\right).
$$

Proposition 4 is a special case of the {\it complementation principle} for the
discrete determinantal processes which is due to S.~Kerov, see \cite{BOO, A.3}.

Observe that Proposition 2 and  Propositions 3 and 4 with similar statements
regarding $\p^{(m)}_v$ and $\q^{(m)}_v$, imply that all the diagonal minors of
the  matrix   $$ K_u^{(m)}=\Vert K_u^{(m)}(x_i,x_j)\Vert_{i,j=0,\dots,M}
$$
are equal to the corresponding diagonal minors of the matrix
$$ 
I-K_v^{(M-m+1)}=\Vert\delta_{ij}- K_v^{(M-m+1)}(x_i,x_j)\Vert_{i,j=0,\dots,M}.
$$
In particular, the diagonal entries of these two matrices are equal.
Looking at 2$\times$2 diagonal minors, we then conclude that
$$
K_u^{(m)}(x,y)=\pm K_v^{(M-m+1)}(x,y)
$$
for all $x\ne y$, $x,y\in X$. (Here we used the fact that both matrices are symmetric.) 

An obvious guess is that the matrices $K_u^{(m)}$ and $I-K_v^{(M-m+1)}$ are conjugate, and the conjugation matrix is diagonal with diagonal entries equal to $\pm 1$. This guess turns out to be correct. 

Set 
$$
D=\operatorname{diag}(\epsilon(x_0),\epsilon(x_1),\dots,\epsilon(x_M)),
$$
where $\epsilon(x)$ was defined in Theorem 1. 

\proclaim{Theorem 5} Under the above notation, for any $m=0,1,\dots, M$,
$$
K_u^{(m)}=D(I-K_v^{(M-m+1)})D,
$$
where the functions $u$ and $v$ satisfy \tht{1}. 
\endproclaim
\demo{Proof}  The equality of the diagonal entries was discussed above: it is
exactly the equality of the first correlation functions of the processes
$\p_{u}^{(m)}$ and $\q_v^{(M-m+1)}$, see Propositions 2, 3, 4. To prove the
equality of the off-diagonal entries we employ the well--known
Christoffel--Darboux formula, see, e.g., \cite{Sz}, which implies that, for
$x\ne y$, $$ \gathered K_u^{(m)}(x,y)=\sqrt{u(x)u(y)}\,\frac{a_{m-1}}{a_m
p_{m-1}}\, \frac{P_m(x)P_{m-1}(y)-P_{m-1}(x)P_m(y)}{x-y}\,,\\
K_v^{(M-m+1)}(x,y)=\sqrt{v(x)v(y)}\,\frac{b_{M-m}}{b_{M-m+1} q_{M-m}}\\ \times
\frac{Q_{M-m+1}(x)Q_{M-m}(y)-Q_{M-m}(x)Q_{M-m+1}(y)}{x-y}\,. \endgathered
$$
Then Theorem 1 immediately implies that $K_u^{(m)}(x,y)=-\epsilon(x)\epsilon(y)K_v^{(M-m+1)}(x,y)$, and the proof is complete.\qed
\enddemo

\head 3. Examples
\endhead

Our main reference for this section is \cite{KS}. We use it for the notation
and data on the classical orthogonal polynomials considered below.

\subhead 3.1. Krawtchouk polynomials \endsubhead
Let $X=\{0,1,\dots,N\}$, and 
$$
u(x)=\binom{N}{x} p^x(1-p)^{N-x}=\frac{N!}{x!(N-x)!}\,p^x(1-p)^{N-x}, \quad
x\in X,\quad 0<p<1. 
$$
The polynomials orthogonal with the weight $u(x)$ are called the {\it
Krawtchouk polynomials}, see \cite{KS, 1.10},
$$
P_n(x)=K_n(x;p,N),\quad n=0,1,\dots,N.
$$
The leading coefficient $a_n$ of $P_n$, the square of the norm $p_n$ of $P_n$,
and the explicit formula for $P_n$ are as follows:
$$
a_n=\binom Nn^{-1}\frac {(-1)^n}{n!p^n},\quad 
p_n=\binom Nn^{-1}\left(\frac {1-p}p\right)^n,\quad
P_n(x)={}_2F_1\left(\matrix -n,\,-x\\ -N\endmatrix\,\Bigr|\frac 1p\right).
$$
Observe that
$$
\prod_{\Sb y=0,\dots,N\\ y\ne x\endSb}(x-y)^2=x!^2(N-x)!^2,\quad
x=0,1,\dots,N. 
$$
Thus, the dual (according to Theorem 1) weight $v(x)$ has the form
$$
v(x)=\left(u(x)x!^2(N-x)!^2\right)^{-1}=\frac 1{N!^2 (p(1-p))^N}\,\binom Nx
(1-p)^x p^{N-x}.
$$
We conclude that $Q_n(x)=\const K_n(x;1-p,N).$ An easy calculation
shows that the normalization of Theorem 1
implies that
$$
\const=(-1)^N(1-p)^N N!, \qquad Q_n(x)=(-1)^N(1-p)^N N! K_n(x;1-p,N).
$$
Clearly, $\epsilon(x)=(-1)^{N-x}$, and the claim of Theorem 1 takes the
form $$
K_n(x;p,N)=(-1)^{x}\left(\frac{1-p}p\right)^{x} K_{N-n}(x;1-p,N),\quad
x=0,\dots,M.  
\tag 2
$$

Of course, this identity can be proved directly using the explicit
formula for the Krawtchouk polynomials above. One just needs to use the
transformation formula  
$$
{}_2F_1\left(\matrix a,\, b\\ c\endmatrix \,\Bigr|\, z\right)=
(1-z)^{-b}{}_2F_1\left(\matrix c-a,\, b\\ c\endmatrix \,\Bigr|\,
\frac{z}{z-1}\right).
$$

\subhead 3.2. Hahn polynomials 
\endsubhead
The computation below was shown to me by T.~Koornwinder. 
Let $X$ be as above, and 
$$
u(x)=\binom{\alpha+x}x\binom{\beta+N-x}{N-x}, \qquad \alpha,\beta>-1\ 
\text{  or  }\  \alpha,\beta<-N.
$$
If $\alpha,\beta>-1$ then $u(x)>0$, if $\alpha,\beta<-N$ then $(-1)^Nu(x)>0$.

The orthogonal polynomials corresponding to this weight are called the {\it
Hahn polynomials}, see \cite{KS, 1.5},
$$
P_n(x)=H_n(x;\alpha,\beta,N),\quad n=0,1,\dots,N.
$$
The data are as follows:
$$
\gathered
a_n=\frac{(n+\alpha+\beta+1)_n}{(\alpha+1)_n (-N)_n}\,,\quad
p_n=\frac{(-1)^n (n+\alpha+\beta+1)_{N+1}(\beta+1)_n
n!}{(2n+\alpha+\beta+1)(\alpha+1)_n (-N)_nN!}\,,\\
P_n(x)={}_3F_2\left(\matrix -n,\, n+\alpha+\beta+1,\,-x\\ \alpha+1,\,
-N\endmatrix\,\Bigr|\,1\right).
\endgathered
$$
The dual weight has the form 
$$
\gathered
v(x)=(u(x)x!^2(N-x)!^2)^{-1}\\=\frac{(-1)^N}{(\alpha+1)_N(\beta+1)_N}
\binom{(-\beta-N-1)+x}x\binom{(-\alpha-N-1)+N-x}{N-x}.
\endgathered
$$
Thus, $Q_n(x)=\const H_n(x;-\beta-N-1,-\alpha-N-1,N).$ Computation of the
normalization constant yields
$$
\const=(-1)^N(\beta+1)^N,\quad Q_n(x)=(-1)^N(\beta+1)^N
H_n(x;-\beta-N-1,-\alpha-N-1,N). $$
The claim of Theorem 1 takes the form
$$
H_n(x;\alpha,\beta,N)=\frac{(-\beta-N)_x}{(\alpha+1)_x}\,
H_{N-n}(x;-\beta-N-1,-\alpha-N-1,N)
\tag 3
$$
for all $x=0,1,\dots,N$.

A direct proof of \tht{3} follows from the transformation formula
$$
{}_3F_2\left(\matrix a,\, b,\,c\\ d,\,
e\endmatrix\,\Bigr|\,1\right)=\frac{\Gamma(d)\Gamma(d+e-a-b-c)}{\Gamma(d+e-a-b)\Gamma(d-c)}\,
{}_3F_2\left(\matrix e-a,\, e-b,\,c\\ d+e-a-b,\,
e\endmatrix\,\Bigr|\,1\right),
$$
see \cite{PBM, 7.4.4(1)}, \cite{Ba, 3.6}.

The limit transition $\alpha=pt$, $\beta=(1-p)t$, $t\to\infty$, see \cite{KS,
2.5.3}, brings \tht{3} to \tht{2}.


\Refs 
\widestnumber\key{AAA}

\ref\key Ba
\by W.~N.~Bailey
\book Generalized hypergeometric series
\publ Cambridge Univ. Press
\publaddr London
\yr 1935
\endref

\ref\key BOO
\by A.~Borodin, A.~Okounkov and G.~Olshanski
\paper Asymptotics of Plancherel measures for symmetric groups
\jour J. Amer. Math. Soc.
\vol 13
\issue 3
\yr 2000
\pages 481--515; {\tt math/9905032}
\endref

\ref\key BO1
\by A.~Borodin and G.~Olshanski
\paper Distributions on partitions, point processes, and the hypergeometric 
kernel
\jour Commun. Math. Phys.
\vol 211
\yr 2000
\pages 335--358; {\tt math/9904010}
\endref

\ref\key BO2
\bysame
\paper $z$--Measures on partitions, Robinson--Schensted--Knuth
correspondence, and $\beta=2$ random matrix ensembles
\paperinfo in Random matrices and their applications. MSRI
Publications Vol. 40, 2001; {\tt math/9905189}
\endref

\ref\key BO3
\bysame
\paper Harmonic analysis on the infinite-dimesional unitary group
\paperinfo In preparation
\endref

\ref\key DVJ
\by D.~J.~Daley, D.~Vere--Jones 
\book An introduction to the theory of point processes 
\bookinfo Springer series in statistics 
\publ Springer 
\yr 1988 
\endref 

\ref\key Dy
\by F.~J.~Dyson
\paper Statistical theory of the energy levels of complex systems I,
II, III 
\jour J. Math. Phys. 
\vol 3
\yr 1962
\pages 140-156, 157-165, 166-175
\endref

\ref\key Ga
\by M.~Gaudin  
\paper Sur la loi limite de l'espacement de valuers propres d'une matrics aleatiore
\jour Nucl. Phys.
\vol 25
\pages 447--458
\yr 1961
\endref

\ref\key GM
\by M.~Gaudin and M.~L.~Mehta
\paper On the density of eigenvalues of a random matrix
\jour Nucl. Phys.
\vol 18
\pages 420--427
\yr 1960
\endref

\ref\key Ha
\by M.~Hall
\book Combinatorial theory
\publ Blaisdell Pub. Co.
\publaddr Waltham, Mass.
\yr 1967
\endref

\ref\key J1
\by K.~Johansson
\paper Shape fluctuations and random matrices
\jour Commun. Math. Phys.
\vol 209
\yr 2000
\pages 437--476;
{\tt math/9903134}
\endref

\ref\key J2
\bysame
\paper Discrete orthogonal polynomial ensembles and the Plancherel
measure 
\paperinfo Preprint, 1999;
{\tt math/9906120}
\endref

\ref\key J3
\bysame
\paper Non-intersecting Paths, Random Tilings and Random Matrices
\paperinfo Preprint, 2000; math/0011250
\endref

\ref\key KS
\by R.~Koekoek and R.~F.~Swarttouw
\paper The Askey--scheme of hypergeometric orthogonal polynomials
and its $q$-analogue \paperinfo available via   
{\tt ftp://ftp.twi.tudelft.nl/\~{}koekoek}
\endref


\ref\key Ma
\by O.~Macchi
\paper The coincidence approach to stochastic point processes
\jour Adv. Appl. Prob.
\vol 7
\yr 1975
\pages 83--122
\endref 

\ref \key Me 
\by M.~L.~Mehta
\book Random matrices
\publ 2nd edition, Academic Press, New York
\yr 1991
\endref

\ref\key NW
\by T.~Nagao, M.~Wadati
\paper Correlation functions of random matrix ensembles related to
classical orthogonal polynomials
\jour  J. Phys. Soc. Japan  \vol 60
\issue 10\yr 1991\pages 3298-3322
\endref  

\ref\key PBM
\by A.~P.~Prudnikov, Yu.~A.~Brychkov, O.~I.~Marichev
\book Integrals and series. Vol. 3: More special functions
\publ Gordon and Breach
\yr 1990
\endref

\ref\key So
\by A.~Soshnikov
\paper Determinantal random point fields.
\paperinfo Russian Math. Surveys, to appear; 
{\tt math/0002099}
\endref

\ref\key Sz
\by G.~Szeg\"o
\book Orthogonal polynomials
\bookinfo AMS Colloquium Publications {\bf XXIII}
\publ Amer. Math. Soc. 
\yr 1959
\publaddr N.Y.
\endref

\endRefs
\enddocument